\documentclass{m2an}

\usepackage{epsfig,graphics,subfigure,psfrag,amsmath,amssymb}\usepackage{amsfonts}
\usepackage{latexsym,wrapfig,picinpar,picins}
\usepackage{indentfirst}
\usepackage{latexsym,bm}
\usepackage{amsthm}
\usepackage{graphicx}
\usepackage{mathrsfs}
\usepackage{epsfig}
\usepackage{fancyref}
\usepackage{dcolumn}
\usepackage{color}
\usepackage{picins}
\usepackage{pstricks}
\usepackage{epstopdf}
\usepackage{subfigure}
\usepackage{multicol}
\usepackage{caption}
\usepackage{geometry}
\usepackage{cases}
\usepackage{amsfonts}
\usepackage{amsopn}
\usepackage{enumerate}
\usepackage{threeparttable}
\usepackage{multirow}
\usepackage{setspace}

\begin{document}
\renewcommand{\theequation}{\arabic{section}.\arabic{equation}}
\setlength\parindent{2em}
\def\open#1{\setbox0=\hbox{$#1$}
\baselineskip = 0pt \vbox{\hbox{\hspace*{0.4 \wd0}\tiny
$\circ$}\hbox{$#1$}} \baselineskip = 10pt\!}
\newtheorem{thm} {Theorem}[section]
\newtheorem{prop}[thm]{Proposition}
\newtheorem{lem}[thm]{Lemma}
\newtheorem{cor}[thm]{Corollary}
\newtheorem{equ}[thm]{Equation}
\newtheorem{defn}[thm]{Definition}
\newtheorem{notation}[thm]{Notation}
\newtheorem{example}[thm]{Example}
\newtheorem{conj}[thm]{Conjecture}
\newtheorem{prob}[thm]{Problem}
\newtheorem{rem}[thm]{Remark}
\title{Numerical Complete Solution for Random Genetic Drift by  Energetic Variational Approach}
%
\author{Chenghua Duan}\address{Department of Mathematics, Soochow University, Suzhou 215006, China; e-mail:  dch3884586@sina.com}
\author{Chun Liu}\address{Department of Mathematics, Pennsylvania State University, University Park, PA, 16802,  USA; e-mail:  liu@psu.edu}
\author{Cheng Wang}\address{Department of Mathematics, University of Massachusetts, Dartmouth, North Dartmouth, MA, 02747-2300, USA; e-mail: cwang1@umassd.edu}
\author{Xingye Yue }\address{Department of Mathematics, Soochow University, Suzhou 215006, China; e-mail:  xyyue@suda.edu.cn}

\date{...}
\begin{abstract}In this paper, we focus on numerical solutions for random genetic drift problem, which is governed by a degenerated convection-dominated parabolic equation. Due to the fixation phenomenon of genes,   Dirac delta singularities will develop at boundary points as time evolves. Based on an energetic variational approach (EnVarA), a balance between the maximal  dissipation principle (MDP) and least action principle (LAP), we obtain the trajectory equation. In turn, a numerical scheme is proposed using a convex splitting technique, with the unique solvability (on a convex set) and the energy decay property (in time) justified at a theoretical level. Numerical examples are presented for cases of pure drift and drift with semi-selection. The remarkable advantage of this method is its ability to catch the Dirac delta singularity close to machine precision over any equidistant grid.\end{abstract}
%
%
\maketitle
\section*{Introduction}

\indent Random genetic drift is the phenomenon that the frequency of a gene variant (allele) in a population changes at the next generation due to random sampling.    The process of random genetic drift  plays an important role in the molecular evolution \cite{J. F. Crow(1970)} and the behavior of genes in a population with a finite size \cite{M. Kimura(1983)}. From the view-point of population genetics, the most elementary step in the evolution is the change of gene frequencies. The notion and technique of random genetic drift have been widely applied to medical science \cite{A. Traulsen(2013)}  and other fields.

\indent We consider a population with a finite size, which  can generally cause the random genetic drift. The  change
in gene frequencies is treated as a stochastic process, which was first introduced by  Fisher  \cite{R. A. Fisher(1922)}.  Under the assumption that generations do not overlap and each copy of gene in the new generation is chosen independently at random from all copies in the old generation, the mathematical model of genetic drift is labeled as the Wright-Fisher Model, introduced by  Fisher \cite{R. A. Fisher(1930)} and Wright \cite{S. Wright(1945)}, and developed by Kimura \cite{M. Kimura(1955a)}. This mathematical model
is a formulation based on a discrete-time Markov chain. The model involves two alleles: $A$ and $a$ in a population with a fixed size $N_{e}$. The quantities $Y_{t}$ and $f$ denote the proportion of $A$ at generation  $t$ in the population and its probability distribution, respectively. Assume that the number of gene $A$ is $n$ at generation $t+1$, which is $m$ at the last generation, the transition probability is given by
$$P\big(Y_{t+1}=\frac{n}{2N_{e}}\Big|{Y}_{t}=\frac{m}{2N_{e}})={2N_{e} \choose n}(\frac{m}{2N_{e}})^{n} (1-\frac{m}{2N_{e}})^{(2N_{e}-n)} , $$
under the circumstance that there is no factor such as mutation, migration and selection and the only evolutionary force is genetic drift. We get the distribution of probability at generation $t+1$ by the Markov chain: $f_{t+1,n}=\sum\limits_{m=1}\limits^{2N_{e}} W_{n,m}f_{t,m}$, where $W_{n,m}$ is transition probability.
 We approximate $Y_{t}$ and $f_{t,n}$ to $x(t)$ and $f(x,t)$, respectively. Kimura \cite{M. Kimura(1955a), M. Kimura(1964), L. Zhao(2013)} showed that for \emph{ pure\ drift} (the only evolutionary force is genetic drift), $f(x,t)$  obeys the diffusion equation:
\begin{equation}\label{equ:diffu}
  \frac{\partial}{\partial t}f(x,t)=\frac{1}{4N_{e}}\frac{\partial^{2}}{\partial x^{2}}(x(1-x)f(x,t)), \ x\in (0,1), \ \ t>0, 
\end{equation}
where $N_{e}$ is the population size. Moreover, if   mutation, migration and selection effects are involved, the model becomes
\begin{equation}\label{equ:sel}
  \frac{\partial}{\partial t}f(x,t)=\frac{1}{4N_{e}}\frac{\partial^{2}}{\partial x^{2}}( x(1-x)f(x,t) )-\frac{\partial}{\partial x}( M(x)f(x,t) ), \ x\in (0,1),\ \ t>0, 
\end{equation}
where $M(x)$ represents the deterministic part of gene frequency dynamics and is typically taken as a polynomial in $x$, whose coefficients depend on mutation rates, migration rates and selection coefficients.

We take the zero current boundary condition
\begin{equation} \label{eqinit} \big\{ \frac{1}{4N_{e}}\partial_{x}[x(1-x)f(x,t)]-M(x)f(x,t)\big\}\mid_{x=0,1}=0, \ \ t>0, $$   with $M(x)=0$ for pure drift and a initial state $$f(x,0)=f_{0}(x)=\delta(x-x_0), \end{equation}
which means that at initial time, the proportion of Gene $A$ is $x_0\in [0,1]$.

  A complete solution, i.e., the total probability is equal to unity at any time, develops sharp spikes (Dirac delta singularities) at the two boundary 0 and 1. When the sharp spikes appear, they signal gene loss or gene fixation: either all copies of Gene $A$ are finally lost,   or all individuals carry $A$ (Gene $a$ is totally lost).
  A complete solution  is essential in Wright Fisher model, because the complete solution can include all possible outcomes whenever fixation and loss are possible, and can be extremely close correspondence with Wright-Fisher model.
  
For the pure drift case, it has been shown  that this system keeps the conservation of the total probability and expectation, and    
$f(x,t) \rightarrow (1-x_0) \delta(x) + x_0 \delta(x-1)$, as $t\rightarrow \infty$ which means that there is a probability of $x_0$ that the  fixation occurs at Gene $A$ and a probability of $1-x_0$ that the   fixation occurs at Gene $a$  \cite{M. Chen(2014), A. J. McKane(2007), T. D. Tran(2013)}.

When  considering  an unlinked locus with two alleles subjects to the semi-dominant selection with strength $s$ ($\mid s\mid\ll 1$),  we take $M(x)=sx(1-x)$  as in \cite{M. Kimura(1955a), M. Kimura(1964)}. 
In this case, the probability of ultimate fixation of Gene $A$ from an initial expectation $x_{0}$  is
$P_{fix}(x_{0})=\frac{1-e^{-4N_{e}sx_{0}}}{1-e^{-4N_{e}s}}$ \cite{M. Kimura(1962), L. Zhao(2013)}.

\indent However, except for a few special cases, we could not get explicit solutions. The numerical approaches are needed to obtain the approximate solutions for the differential equation. Some attempts have been made by Kimura \cite{M. Kimura(1955a)}, Barakat and Wagener \cite{R. Barakat(1978)}
 and Wang \cite{Y. Wang(2004)},
while the total probability is smaller than unity and it was also a hard work to simulate the general case  including natural selection, mutation and migration. Zhao  et\ al. \cite{L. Zhao(2013)} obtained a complete numerical solution by finite volume method (FVM) for a neutral locus and semi-selection. 
In \cite{M. Chen(2014)}, Xu et al. discussed three classical numerical schemes  which are stable but lead to different steady state solutions. Only one of  the schemes gives a true complete numerical solution and any scheme with numerical viscosity should be avoided. Therefore, a very careful analysis for the numerical scheme is necessary.


\indent{In} this paper, we propose a new scheme based on energetic variational approach (EnVarA). Combining the least action principle (LAP) and maximal  dissipation principle (MDP), we first obtain the trajectory equation for the Wright-Fisher model. In turn, a convex-splitting technique is applied to construct a numerical scheme that is unique solvable on a convex domain and keeps the property of energy decay in time. The numerical scheme can assure the conservation of the  total probability, i.e., a complete solution is obtained.  Numerical examples demonstrate that we can get a complete solution and true probability of fixation. In comparison with the FVM schemes in \cite{M. Chen(2014), L. Zhao(2013)}, the new method has a significant advantage on the approximation to the delta singularity. Over an equidistant mesh with step size $h$, standard finite difference methods or FVMs only present an approximation of scale $O(1/h)$ to  delta singularity, while the scheme here may give an approximation of scale $O(1/\varepsilon)$ with small positive $\varepsilon$ close to  the machine precision.

    The paper is organized as follows. The details of EnVarA for Wright-Fisher model are shown in Section \ref{sec:2}. In Section \ref{sec:3}, the numerical scheme is constructed. Then numerical examples  are presented in Section \ref{sec:4}.

\section{Variational approach for the Wright Fisher model}\label{sec:2}
\setcounter{equation}{0}
The primary goal of this section is to derive the constitutive relation of the Wright Fisher model. We first introduce  EnVarA briefly. The original work was given by Onsager \cite{L. Onsager(1931)}, and then it was improved by Rayleigh \cite{J. W. Strutt(1873)}. This method has been applied to many physical and biological problems in recent years, for instance \cite{B. Eisenberg(2010),  Q. Du(2009), X. F. Yang(2006)}. 
In the Wright-Fisher model, $x\in[0,1]$ and $f(x,t)\geq0$ can be viewed as the position of particles and the density of $x$ at time $t$, respectively. We first introduce  the different coordinate systems.
 \begin{defn}
 Suppose that $\Omega_{0}^{X}$, $\Omega_{t}^{x}$ $\subset\mathbb{R}^{m}$, $m\in\mathbb{N}^{+}$, are domains with smooth boundary and time $t>0$, and $\textbf{u}=(u_{1},...,u_{m})$ is a smooth vector field in $\mathbb{R}^{m}$. The flow map $x(X,t):\Omega_{0}^{X}\rightarrow\Omega_{t}^{x}$ is defined as a solution of:
 \begin{equation}\label{equ:v}
 \left\{
 \begin{aligned}
   & \frac{d}{dt}x(X,t)=\textbf{u}(x(X,t),t),\ \ t>0 , \\
   & x(X,0)=X ,
 \end{aligned}
 \right.\
 \end{equation}
 where $X=(X_{1},...,X_{m})\in\Omega_{0}^{X}$ and $x=(x_{1},...,x_{m})\in\Omega_{t}^{x}$. In turn, the coordinate system $X$ is called the Lagrangian coordinate and the coordinate system $x$ is called Eulerian coordinate.
 \end{defn}
EnVarA is obtained by the combination of the statistical physics and nonlinear thermodynamics. First, we define total energy $$E^{total}:=\mathcal{K}+\mathcal{H},$$ where $\mathcal{K}$ is the kinetic energy and $$\mathcal{H}:=\mathcal{U}-T\mathcal{S}$$ is the Helmholtz free energy containing the internal energy $\mathcal{U}$, temperature $T$
 and entropy $S$.
In  an isothermal system without external force, the total energy dissipation law holds:
 $$\frac{d}{dt}E^{total}=-\Delta,$$
 where $\Delta\geq0$ is the entropy product. 

Subsequently, the least action principle (LAP) is applied: the trajectory of particles $X$ from $x(X,0)$ at time $t=0$ to
$x(X,t^{*})$ at a given time $t^{*}$ in a Hamiltonian system are those which minimize the action functional defined by
$$\mathcal{A}(x(X,t)):=\int^{t^{*}}_{0}\mathcal{L}(x(X,t),x_{t}(X,t))dt,$$
where $\mathcal{L}:=\mathcal{K}-\mathcal{H}$ is the Lagrangian functional of a conservative system and $x(X,t)\in\Omega_{t}^{x}$, $t>0$.
 Moreover, in a non-Hamiltonian system here, taking variational of the action functional with respect to $x$, we get the conservation force
 $$F_{con}=\frac{\delta\mathcal{A}}{\delta x}.$$
%
Next, we treat the dissipation part with maximum dissipation principle (MDP).
Taking variational of $\Delta$ with respect to the velocity $\textbf{u}$ involved in (\ref{equ:v}), we have the dissipative force
$$F_{dis}=\frac{\delta\frac{1}{2}\Delta}{\delta\textbf{u}},$$
where the factor $\frac{1}{2}$ comes from a linear reponse assumption, i.e., $\Delta$ is quadratic  function  of $\textbf{u}$ and  $F_{dis}$ is linear in $\textbf{u}$ \cite{R. Kubo(1976)}. According to the Newton's force balance law: $$F_{con}=F_{dis},$$
we obtain constitutive relation. Onsager's approach \cite{L. Onsager(1931), L. Onsager1(1931)} is the key point for such conclusions.

Now we revisit  the  Wright-Fisher model with a positive initial state in a context of EnVarA. By rescaling the time,   (\ref{equ:diffu}) \eqref{eqinit} becomes:
\begin{eqnarray}
 &&    \partial_t f+ \partial_x(f\textbf{u})=0 ,  \label{eqcm}\\
 && f\textbf{u} =-  \partial_x \big(x(1-x)f\big) , \label{eqwf} \\
   && f(x,0)=f_{0}(x)>0,\ x\in [0,1] ,  \\
     && \partial_{x}(x(1-x)f)\mid_{x=0,1}=0, \ t>0.  \label{eqbc}
    \end{eqnarray}

\begin{lem} \label{lem-dissip}
  $f(x,t)$ is the solution of \eqref{eqcm}-\eqref{eqbc} if and only if
  $f$ satisfies the corresponding energy dissipation law
  \begin{equation}\label{equ:energylaw}
   \frac{d}{dt}\int_{0}^{1}f \ln(x(1-x)f)dx=-\int_{0}^{1}\frac{f}{x(1-x)}|\textbf{u}|^{2}dx.
  \end{equation}
    \end{lem}

\noindent\textbf{Proof}:
 We first prove that the energy dissipation law (\ref{equ:energylaw}) holds if $f$ is the solution of  \eqref{eqcm}-\eqref{eqbc}.
Multiplying  by $1+\ln{(x(1-x)f)}$  and integrating on both sides of (\ref{eqcm}), we get
$$\int_{0}^{1}\Big(1+\ln\big(x(1-x)f\big)\Big)\partial_t f dx=\int_{0}^{1}\Big(1+\ln\big(x(1-x)f\big)\Big)\partial_{xx}(x(1-x)f)dx .$$
By integration by parts, we have
\begin{subequations}
  \begin{align}
  \frac{d}{dt}\int_{0}^{1} f\ln\big(x(1-x)f\big) dx
=&-\int_{0}^{1}\frac{\partial}{\partial x}(x(1-x)f)\frac{\frac{\partial}{\partial x}(x(1-x)f)}{x(1-x)f}dx \nonumber\\
  =&-\int_{0}^{1}\frac{f}{x(1-x)}|\textbf{u}|^{2}dx.  \nonumber
  \end{align}
\end{subequations}

Next we  can derive  \eqref{eqwf}  from  the   energy dissipation law  (\ref{equ:energylaw}) by EnVarA, while (\ref{eqcm}) is the conservation law which is assumed to be true.

Note that in  Lagrangian coordinate, there exists an explicit formula for the solution of the conservation law (\ref{eqcm}),
 \begin{equation}\label{equ:conservationL}
  f(x(X,t),t)=\frac{f_{0}(X)}{\frac{\partial x(X,t)}{\partial X}},
\end{equation}
where $f_{0}(X)$ is the  initial function and $\frac{\partial x(X,t)}{\partial X}$ is \emph{deformation} \emph{gradient}, which is the Jacobian matrix of the map: $X\rightarrow x(X,t)$.

\begin{itemize}
\item
  The {\bf total energy} of the Wright-Fisher model is given by
  \begin{equation} \label{tot-energ-e} E^{total}=\mathcal{H}=\int_{0}^{1}f\ln(x(1-x)f) dx. \end{equation}
\item LAP step. \ \
With (\ref{equ:conservationL}), the action functional in Lagrangian coordinate becomes
$$\mathcal{A}(x)=\int^{t^{*}}_{0}(-\mathcal{H}) dt = -\int^{t^{*}}_{0}\int_{0}^{1}f_{0}(X)\ln\left(x(1-x)\frac{f_{0}(X)}{\frac{\partial x(X,t)}{\partial X}}\right)dXdt , $$
where $t^{*}>0$ is a given  terminal time.
Thus for any test function $y(X,t)=\widetilde{y}(x(X,t),t)\in C_{0}^{\infty}((0,1)\times (0,t^*))$ and $\epsilon\in \mathbb{R}$, taking the variational of $\mathcal{A}(x)$ with respect to $x$, we get
\begin{subequations}
  \begin{align}
  \frac{d}{d\epsilon}\bigg|_{\epsilon=0}\mathcal{A}(x+\epsilon y)%
  &=-\int^{t^{*}}_{0}\int_{0}^{1}\left(f_{0}(X)\frac{1-2x}{x(1-x)}+\frac{\partial}{\partial X}\left(\frac{f_{0}(X)}{\frac{\partial x}{\partial X}}\right)\right)y dXdt \nonumber \\
  &=-\int^{t^{*}}_{0}\int_{0}^{1}\left(f\frac{1-2x}{x(1-x)}+\frac{\partial f}{\partial x}\right)\widetilde{y}dxdt . \nonumber
  \end{align}
\end{subequations}
  Then we obtain the conservation force
  $$F_{con}= \frac{\delta\mathcal{A}}{\delta x} =-\left(f\frac{1-2x}{x(1-x)}+\frac{\partial f}{\partial x}\right)=-\frac{1}{x(1-x)}\frac{\partial}{\partial x}\big(x(1-x)f\big), $$ in Eulerian coordinate, and
   $$F_{con}=-\left(f_{0}(X)\frac{1-2x}{x(1-x)}+\frac{\partial}{\partial X}\left(\frac{f_{0}(X)}{\frac{\partial x}{\partial X}}\right)\right),$$ in Lagrangian coordinate.

\item  MDP step.  \ \   Let the entropy production
$\Delta=\int_{0}^{1}\frac{f}{x(1-x)}|\textbf{u}|^{2}dx$. Taking the variational of $\frac 12 \Delta$ with respect to $\textbf{u}$, we have the dissipation force
$$F_{dis}=\frac{\delta\frac{1}{2}\Delta}{\delta\textbf{u}}=\frac{f}{x(1-x)}\textbf{u}, $$ in Eulerian coordinate, and
$$F_{dis}=\frac{\delta\frac{1}{2}\Delta}{\delta x_{t}}=\frac{f_{0}(X)}{x(1-x)}x_{t}, $$ in Lagrangian coordinate.\\
\item Force  balance step. \ \ We have, in Lagrangian coordinate, that
\begin{equation} \label{trajLa}
\frac{f_{0}(X)}{x(1-x)}x_{t}=-\frac{\partial}{\partial X}\left(\frac{f_{0}(X)}{\frac{\partial x}{\partial X}}\right)-f_{0}(X)\frac{1-2x}{x(1-x)},
\end{equation}
and in Eularian coordinate, we have
\begin{equation}\label{ne}
    \begin{aligned}
    \frac{f(x,t)}{x(1-x)}\textbf{u} =  -\frac{1}{x(1-x)}\frac{\partial}{\partial x}\Big( x(1-x)f(x,t)\Big),
    \end{aligned}
\end{equation}
which is exactly \eqref{eqwf}.  $\hfill\Box$\\
\end{itemize}

\begin{rem}  \label{reminit} There is an assumption that the initial state is positive in the above lemma. Otherwise, if $f_0(X)=0$ for some $X\in (0,1)$, the argument  above would be not valid any more. For example, in \eqref{trajLa}, the velocity $x_t$ could be indefinite for points such that $f_0(X)=0$. Note that in the real model, the initial state \eqref{eqinit} is $f_0= \delta(x-x_0)$, almost  zero everywhere. To deal with this case, we consider two models with positive initial states $f_{0,1}, f_{0,2}$ such that $f_0= f_{0,1}-f_{0,2}$ and correspondingly we have $f = f_1(x,t)-f_2(x,t)$.
\end{rem}

\begin{rem} \label{remalg0} What we really get by EnVarA is  \eqref{trajLa}, which contains all the physics involved in this model.  If we can solve \eqref{trajLa} to get the trajectory $x(X,t)$, substituting it into \eqref{equ:conservationL}, we obtain the   solution $f(x,t)$ to \eqref{eqcm}-\eqref{eqbc}. So in the following sections, we focus on numerical solution to \eqref{trajLa}.
\end{rem}

To this purpose, we should first settle the initial and boundary condition for \eqref{trajLa}. From \eqref{eqbc} and \eqref{eqwf}, we have $x_t(0,t)=x_t(1,t)=0$, for $t>0$. That means that a Dirichlet boundary condition should be subject to as $x(0,t)=0, x(1,t)=1$, for $t>0$. So the trajectory problem is
\begin{equation}\label{equ:crp}
   \left\{
    \begin{aligned}
     &\frac{f_{0}(X)}{x(1-x)}\partial_t x=-\frac{\partial}{\partial X}\left(\frac{f_{0}(X)}{\frac{\partial x}{\partial X}}\right)-f_{0}(X)\frac{1-2x}{x(1-x)},\ X\in (0,1),\ t>0, \\
     & x(X,0)=X,\ X\in [0,1], \\
     & x(0,t)=0,\ x(1,t)=1,\ t>0.
    \end{aligned}
    \right. \
\end{equation}

\section{Numerical methods for trajectory equation}
\setcounter{equation}{0}
\label{sec:3}
 In this section, we consider numerical methods  for (\ref{equ:crp}). 
\subsection{ A semi-discrete  scheme in time and optimal transport}
System (\ref{equ:crp}) can be viewed as a gradient flow associated with the total energy of
 \begin{equation} \label{eqtotEn} E^{total}=\int_{0}^{1}f_{0}(X)\ln\Big(\frac{f_{0}(X)}{\frac{\partial x}{\partial X}}\Big)dX+\int_{0}^{1}f_{0}(X)\ln\big(x(1-x)\big)dX,
\end{equation}
which is just the counterpart in Lagrangian coordinate of total energy \eqref{tot-energ-e} of the system \eqref{eqcm}-\eqref{eqbc} and can be split into  convex and concave parts, that is $E^{total}=E_{c}-E_{e}$, where both $E_{c}$ and $E_{e}$ are convex. The canonical splitting is  $E_{c} =\int_{0}^{1}f_{0}(X)\ln\Big(\frac{f_{0}(X)}{\frac{\partial x}{\partial X}}\Big)dX$ and $E_{e} =-\int_{0}^{1}f_{0}(X)\ln\big(x(1-x)\big)dX$.
 The convex splitting was first exploited by D. J. Eyre in \cite{{D. J. Eyre(1998)}} to craft energy stable numerical schemes for the Allen-Cahn and Cahn-Hilliard equations. The basic idea is to treat the convex part implicitly while to treat the concave part explicitly. Then a {\bf  semi-discrete  scheme} for (\ref{equ:crp}) is proposed as follows
\begin{equation}\label{equ:num}
 \frac{f_{0}(X)}{x^{n}(1-x^{n})}\frac{x^{n+1}-x^{n}}{\tau}=
-\frac{\partial}{\partial X}\left(\frac{f_{0}(X)}{\frac{\partial x^{n+1}}{\partial X}}\right)-f_{0}(X)\frac{1-2x^{n}}{x^{n}(1-x^{n})} ,
\end{equation}
where $\tau$ is the time step and $x^{n}=x(X,t^{n})$ is the solution at time $t^{n}=n\tau$, $n\in\mathbb{N}^{+}$.

\begin{rem}
   \eqref{equ:num} is also a Variational Particle Scheme. We explain the fact in the framework of optimal transport theory. Let $\Omega=[0,1]$. We denote by $\mathscr{P}(\Omega)$ the space of $\mathfrak{L}^{1}$ measure on $\Omega$, non-negative functions with unit integral and finite second moments, where $\mathfrak{L}^{1}$ is the Lebesgue measure. $f^{n}\in\mathscr{P}(\Omega)$ is the approximation to solution of equation \eqref{eqcm}-\eqref{eqwf} at time $t_{n} = n\tau$, $n\in\mathbb{N}$. We fix a reference density $f^{0}$ and consider a time-dependent family of transport maps $x(\cdot ,t^{n}):[0,1]\rightarrow [0,1]$ such that $x(\cdot,t^{n})\#f^{0}=f ^{n}\mathfrak{L}^{1}$ for all $n\in\mathbb{N}^{+}$, where $\#$ denotes the push-forward of measures.

Then the map from $x^n$ to $x^{n+1}$ is an optimal transport in the sense that $x^{n+1}$ is the minimizer of
  the cost functional
$$F(x) :=\int_{0}^{1}\frac{1}{2\tau}\frac{f_{0}(X)}{x^{n}(1-x^{n})}|x-x^{n}|^{2}+
f_{0}(X)\ln\left(\frac{f_{0}(X)}{\frac{\partial x}{\partial X}}\right)+ f_{0}(X)\frac{1-2x^{n}}{x^{n}(1-x^{n})}xdX.$$
Some relevant descriptions on optimal transport can be found in \cite{M. Westdickenberg(2010)}.
\end{rem}

\subsection{The fully discrete scheme}
We begin with the definition of inner-product, difference operators and summation-by-parts in one dimension.
Let $h=\frac{1}{N}$, $N\in\mathbb{N}^{+}$ be the spatial step. Denote by $X_{r}=X(r)=r h$, where $r$ takes on integer and half integer values.  Let $\mathcal{E}_{N}$ and $\mathcal{C}_{N}$ be the spaces of functions whose domains are $\{X_{i}\ |\ i=0,...,N\}$ and $\{X_{i-\frac{1}{2}}\ |\ i=1,...,N\}$ respectively. In component form, these functions are identified via
$l_{i}=l(X_{i})$, $i=0,...,N$, for $l\in\mathcal{E}_{N}$, and $\phi_{i-\frac{1}{2}}=\phi(X_{i-\frac{1}{2}})$,  $i=1,...,N$, for $\phi\in\mathcal{C}_{N}$.

\indent{Let} $l$, $g\in\mathcal{E}_{N}$ and $\phi$, $\psi\in\mathcal{C}_{N}$. We define the ``inner-product" on space $\mathcal{E}_{N}$ and $\mathcal{C}_{N}$  respectively as
\begin{equation}\label{def:inner-integer}
[l\big|g]=h \sum\limits_{i=1}^{N-1} l_{i} g_{i},
\end{equation}
\begin{equation}\label{def:inner-half} (\phi\big|\psi)=h\sum\limits_{i=1}^{N}\phi_{i-\frac{1}{2}}\psi_{i-\frac{1}{2}} .
\end{equation}
\indent{The} difference operator $D_{h}:\mathcal{E}_{N}\rightarrow\mathcal{C}_{N}$ and $d_{h}:\mathcal{C}_{N}\rightarrow\mathcal{E}_{N}$, and the average operator $A:\mathcal{E}_{N}\rightarrow\mathcal{C}_{N}$  can be defined as respectively as
\begin{align}\label{equ:dif1}
& (D_{h}l)_{i-\frac{1}{2}}= (l_{i}-l_{i-1})/h,\ i=1,...,N, \\
&\label{equ:dif2}  (d_{h}\phi)_{i}= (\phi_{i+\frac{1}{2}}-\phi_{i-\frac{1}{2}})/h,\ i=1,...,N-1, \\
&\label{equ:ave1} (Al)_{i-\frac{1}{2}}= (l_{i}+l_{i-1})/2,\ i=1,...,N.
\end{align}
Then we have the following result of summation-by-parts.
\begin{lem}
Let $\phi\in\mathcal{C}_{{N}}$ and $l\in\mathcal{E}_{N}$. Then
$(D_{h}l\big|\phi)=-[l \big|d_{h}\phi]+l_{N}\phi_{N-\frac{1}{2}}-l_{0}\phi_{\frac{1}{2}}$.
\label{lemmasum}
\end{lem}
Let $\textbf{Q}:=\{l \in\mathcal{E}_{N}\ |\ l_{i-1}<l_{i},\ 1\leq i\leq N;\ l_{0}=0,\ l_{N}=1\}$ and its boundary set $\partial\textbf{Q}:=\{l \in\mathcal{E}_{N}\ |\ l_{i-1}\leq l_i,\ 1\leq i\leq N,\ and\  l_{i}=l_{i-1},\ for\ some\ 1\leq i\leq N; l_{0}=0,\ l_{N}=1\}$. Then $\bar{\textbf{Q}}:=\textbf{Q}\cup\partial\textbf{Q}$ is a closed convex set.

The {\bf  fully discrete scheme} is formulated as follows: Given $x^{n} \in\textbf{Q}$, find $x^{n+1}=(x^{n+1}_{0},...,x^{n+1}_{N})\in\textbf{Q}$ such that
\begin{equation}\label{equ:numnum}
 \frac{f_{0}(X_{i})}{x^{n}_{i}(1-x^{n}_{i})}\frac{x^{n+1}_{i}-x^{n}_{i}}{\tau}=
-d_{h}\left(\frac{Af_{0}(X)}{D_{h}x^{n+1}}\right)_{i}
-f_{0}(X_{i})\frac{1-2x^{n}_{i}}{x^{n}_{i}(1-x^{n}_{i})},\ 1\leq i\leq N-1.
\end{equation}
 (\ref{equ:numnum}) is still a nonlinear system. Newton's iteration method can be applied to solve it.

\noindent{\bf Damped Newton's iteration.} \ \ Set $x^{n+1, 0}= x^n$. For $k=0,1,2,\cdots, x^{n+1,k+1}=x^{n+1,k} + \omega(\lambda)\delta_x$ such that
\begin{align}\label{equ:numnum1}
\frac{f_{0}(X_{i})}{x^{n}_{i}(1-x^{n}_{i})}\frac{ \delta_{x_{i}}}{\tau} &- d_h \left(\frac{Af_{0}(X)}{(D_{h}x^{n+1,k})^2} D_h \delta_{x_i}\right)_{i}= -
 \frac{f_{0}(X_{i})}{x^{n}_{i}(1-x^{n}_{i})}\frac{x^{n+1,k}_{i} -x^{n}_{i}}{\tau} \notag \\
&
-d_{h}\left(\frac{Af_{0}(X)}{D_{h}x^{n+1,k}}\right)_{i}
-f_{0}(X_{i})\frac{1-2x^{n}_{i}}{x^{n}_{i}(1-x^{n}_{i})},\ \ 1\leq i\leq N-1.
\end{align}

and
\begin{equation}\label{Newton}
\omega(\lambda)=\left\{
\begin{array}{lcl}
\frac{1}{\lambda} &&{\lambda > \lambda'}\\
\frac{1-\lambda}{\lambda(3-\lambda)} &&{\lambda'\ge\lambda\ge\lambda^{*}}\\
1 && {\lambda < \lambda^{*}},
\end{array}\right.
\end{equation}
where $\lambda^*=2-3^{\frac{1}{2}}$, $\lambda'\in[\lambda^*, 1)$ and
$\lambda(J,x^{n+1,k})=\big(\frac{1}{a}(J'(x^{n+1,k}))^T[J''(x^{n+1,k})]^{-1}J'(x^{n+1,k})\big)^{\frac{1}{2}}$ with $J$  defined in \eqref{DisEnergy}
and $a=h\min\limits_{i}(f_0(X_i))$.

After solving \eqref{equ:numnum},  we finally get the numerical distribution $f(x^{n+1},t^{n+1})$ from \eqref{equ:conservationL} as
\begin{align} \label{numdist}
& f_{i}^{n+1}= \frac{ f_{0}(X_{i}) } {(x_{i+1}^{n+1}-x_{i-1}^{n+1})/(2h)},\ 1\leq i\leq N-1, \ \mbox{and}\\
& f_{0}^{n+1}=\frac{f_{0}(X_{0})} {(x_{1}^{n+1}-x_{0}^{n+1})/h}, \ f_{N}^{n+1}=\frac {f_{0}(X_{N})}{(x_{N}^{n+1}-x_{N-1}^{n+1})/h}. \label{numdist2}
\end{align}

\begin{lem} \label{conservation-mass} The density function $f^{n+1}$ obtained from  \eqref{numdist}-\eqref{numdist2} keeps the conservation law of mass.
\end{lem}
In fact, if we define the initial mass carried by each particle $x_i^0=X_i$ as
\begin{equation} \label{init-mass}
m_i^0= h f_0(X_i), \ 1<i<N; \ \ m_0^0 = \frac h2 f_0(X_0); \ \ m_N^0 = \frac h2 f_0(X_N), \end{equation}
and define the mass carried by particle $x_i^n$ as
 \begin{equation}\label{mass}
m_i^n= \frac{x_{i+1}^n-x_{i-1}^n} 2  f_i^n, \ 1<i<N; \ \ m_0^n = \frac {x_1^n-x_0^n}2 f_0^n; \ \ m_N^n = \frac {x_N^n-x_{N-1}^n}2 f_N^n, \end{equation}
then we readily have from  \eqref{numdist}-\eqref{numdist2} that
$$m_i^n \equiv m_i^0, \ \ 0\leq i\leq N, \ n=1,2,\cdots.$$

\begin{rem} \label{rem:numdelta}
 $x_i(t)=x(X_i,t), 0<i<N$ are the trajectories starting from the particles $X_i$ at time $t=0$. From the governing equation \eqref{equ:crp} or \eqref{equ:numnum}, the motion of these particles is primarily determined by the second term on the right hand side  since this term tends to infinity when the particle approaches to the end points $x=0,1$. In particular, this term tends to negative infinity around the left end $x=0$, while the limit becomes positive infinity around the right end $x=1$.  Therefore, $x_1(t)$ and $x_{N-1} (t)$ will be closer and closer to $x_0(t)\equiv 0$ and $x_N(t)\equiv 1$, respectively.

 Governed  by the continuous model \eqref{equ:crp},  the particles may touch the end points, which means that the Dirac delta singularity occurs for $f(x,t)$ from \eqref{equ:conservationL}. For the discrete model \eqref{equ:numnum}, we find solution $x^{n+1} \in \textbf{Q}$, where $x_i< x_{i+1}$ for $0\leq i< N$. As a result, theoretically $x_1$ and $x_{N-1}$ would never touch the ends. However, in the practical computations, when $x_1^n$ and  $x_0^n=0$ are too close to distinguish from each other under the machine precision, they are bundled up and will be regarded as {\bf one particle} which carries the mass from the original two and will be fixed at the boundary.  This is the signal that the numerical Dirac delta (i.e., the {\bf fixation})  happens. In comparison with the FVMs  in \cite{M. Chen(2014)}, we can now approximate the delta singularity to the scale of $1/\varepsilon$, with $\varepsilon$ close to the machine precision, while by the standard FVMs on equidistance mesh, one can only approximate the delta singularity to the scale of $1/h$ (with the spatial mesh size $h$).
\end{rem}

\noindent{\bf Criteria for particles meet the boundary.}\ \  Though we can choose the machine precision as a criterion to judge whether two particles touch each other, it is not practical. For example, in \eqref{numdist2}, when $x_1^{n+1}-x_0^{n+1}$ is close to machine precision, we will lose all the accuracy of $f_0^{n+1}$. So we will choose a   criterion  with $\varepsilon_0=10^{-10}$ in double precision system  as: \\
 \begin{equation} \label{criteria} \mbox {Criteria:} \left\{ \begin{array}{l}
  \mbox{If\ }  x_i^{n+1}\in \textbf{B}_{l}=[0,\varepsilon_0],  \mbox{it\ will\ be\ fixed\ at\ left\ boundary\ for\ ever,}\\
  \mbox{If\ } x_i^{n+1}\in \textbf{B}_{r}=[1-\varepsilon_0,1],   \mbox{it\ will\ be\ fixed\ at\ right\ boundary\ for\ ever.}
  \end{array}
  \right.
  \end{equation}
  Equivalently, we have a rearrangement on the position of the particles as
  \begin{equation}\label{re-arrange}
 x_{i}^{n+1}=
    \begin{cases}
     0 , & \mbox{if $x_{i}^{n+1}\in\textbf{B}_{l}$} , \\
     1 , & \mbox{if $x_{i}^{n+1}\in\textbf{B}_{r}$} , \\
     x_{i}^{n+1} , & \mbox{  otherwise  .}
    \end{cases}
\end{equation}
  At the next time step, we only need to determine the position of particles from $(\varepsilon_0, 1- \varepsilon_0)$.

With the above rearrangement, the formulas \eqref{numdist2} for the  density function at the boundary points don't work any more. To define the revised formulas, we need to count the total number of particles accumulated at the boundary points. Let
\begin{equation} \label{num-accum}  \begin{cases}
 \mbox{starting\ point}\ & i_s^{n+1} = \max\{\ i\ |x_{i}^{n+1}\in\textbf{B}_{l},\ 0\leq i<N \} , \\
 \mbox{ending\ point}\  &
i_e^{n+1} =\min\{\ i\ |x_{i}^{n+1}\in\textbf{B}_{r},\ 0<i\leq N\}.
\end{cases}
\end{equation}
If $i_s^{n+1} >0$ or $i_e^{n+1} <N$,  there must be some particles which touched the boundary points at time $t^{n+1}$.
Then the revised formula for the density function   $f^{n+1} = (f^{n+1}_{i_s^{n+1}}, f^{n+1}_{i_s^{n+1}+1}, \cdots, f^{n+1}_{i_e^{n+1}})$ become
\begin{align}
& f_{i}^{n+1}= \frac{ f_{0}(X_{i}) } {(x_{i+1}^{n+1}-x_{i-1}^{n+1})/(2h)} = \frac{m^0_i } {(x_{i+1}^{n+1}-x_{i-1}^{n+1})/2},\ {i_s^{n+1}}< i< {i_e^{n+1}},  \label{numdist3} \\
& f_{i}^{n+1}=\frac 2{\varepsilon_0} \sum_{k=0}^{i -1} m^0_k+ \frac{m^0_{i }}{(x_{i +1}-x_{i})/2},\ \mbox{for\ } i=i_s^{n+1}, \mbox{and} \label{numdist4}\\
& f_i^{n+1}=\frac 2{\varepsilon_0}
\sum_{k=i+1}^N m_{k}^0 + \frac{m^0_{i }}{(x_{i}-x_{i-1})/2},\ \mbox{for\ } i= i_e^{n+1}.\label{numdist5}
\end{align}

\begin{rem} The treatment in \eqref{numdist3}-\eqref{numdist5} keeps the conservation law of total mass naturally and  means that only the last fixed particle can feel the free nearest particle inside and the effect of all former fixed particles is confined to the $\varepsilon_0$ neighbor of boundary points.
\end{rem}

Combining all the discussions above together, we can now present the final algorithm as follows.\\

\noindent{\bf Algorithm 2.1. } {\em
\begin{itemize}
\item Initialization. \\ For $0\leq i\leq N$, we get the initial particle position $x_i^0 = X_i$, the initial density distribution $f^0_i = f_0(X_i)$, and the initial mass $m^0_i$ by \eqref{init-mass}.\\
     Set starting point $i_s = 0$ and ending point $i_e = N$.

\item Time Stepping. \\ For $n=0, 1, 2, \cdots $, find the density distribution at next time step $f^{n+1} = (f^{n+1}_{i_s }, f^{n+1}_{i_s +1}, \cdots, f^{n+1}_{i_e})$ by the following procedures.
    \begin{enumerate}
        \item Obtain the position of particles $x^{n+1}_i$, $i_s\leq i \leq i_e$, via solving the fully discrete system \eqref{equ:numnum} by Newton's iteration \eqref{equ:numnum1}, with $x^{n+1}_{i_s} = 0, \ x^{n+1}_{i_e}=1$.
        \item Check whether a particle meets the boundary by the criteria \eqref{criteria}, re-arrange the position by \eqref{re-arrange} and update the starting point $i_s$ and the ending point $i_e$ by \eqref{num-accum} if necessary.
        \item Obtain the density distribution $f^{n+1}$ by \eqref{numdist3}-\eqref{numdist5}.
    \end{enumerate}
\end{itemize}
}

\subsection{Unique solvability  and energy decay of fully discrete scheme}
In this  subsection, we provide some analyses on the unique solvability and energy decay of the fully discrete scheme \eqref{equ:numnum}, and the convergence of  the Newton method \eqref{equ:numnum1} with \eqref{Newton}.

\begin{thm}
\label{lem:unique}
The numerical scheme \eqref{equ:numnum} is unique solvable in  $\textbf{Q}$.
\end{thm}

\noindent\textbf{Proof:} We first consider the following optimization problem: 
\begin{equation}\label{DisEnergy}
\min\limits_{y\in\bar{\textbf{Q}}}\Big\{J(y):=\frac{1}{2\tau}\Big[\frac{f_{0}(X)}
{x^{n}(1-x^{n})}(y-x^{n})\Big|(y-x^{n})\Big]+\Big
(Af_{0}(X)\Big|\ln\big(\frac{Af_{0}(X)}{D_{h}y}\big)\Big)+\Big[ f_{0}(X)\frac{1-2x^{n}}{x^{n}(1-x^{n})}\Big|y\Big]\Big\},
\end{equation}
where $f_{0}(X)\in\mathcal{E}_{N}$ is the initial  distribution  and $x^{n}\in\textbf{Q}$ is the known  position of particles at time $t^{n}$. It is easy to verify that $J(y)$ is a convex function on the closed convex set $\bar{\textbf{Q}}$. Hence there exists a unique minimizer  $x\in\bar{\textbf{Q}}$. We must have  the minimizer  $x\in\textbf{Q}$ since if $y\in\partial\textbf{Q}$, then there exists some $i>0$ such that $(D_h y)_{i-1/2} = (y_i-y_{i-1})/h=0$, and  $J(y)=+\infty$ .

We  first claim  that $x\in\textbf{Q}$ is the minimizer of $J(y)$  if and only if it is a solution of   scheme (\ref{equ:numnum}). Hence the fully discrete scheme \eqref{equ:numnum} has a unique solution.

In fact, if $x\in\textbf{Q}$ is the minimizer of $J(y)$, then for $\forall y\in\bar{\textbf{Q}}$, there exists a sufficiently small $\epsilon_0>0$, such that for any $\epsilon\in(-\epsilon_0, \epsilon_0)$, $x+\epsilon (y-x) \in \textbf{Q}$ since $\textbf{Q}$ is a open set. Then $j(\epsilon) = J(x+\epsilon (y-x))$ achieves its minimal at $\epsilon=0$.  So we have $j'(0)=0$ and  using summation by parts, we obtain

\begin{equation*}
    \frac{1}{\tau}\Big[\frac{f_{0}(X)}{x^{n}(1-x^{n})}(x-x^{n})\Big|y-x\Big]
    +\Big[d_{h}\big(\frac{Af_{0}(X)}{D_{h}x}\big)\Big|y-x\Big]
    +\Big[f_{0}(X)\frac{1-2x^{n}}{x^{n}(1-x^{n})}\Big|y-x\Big]=0,
\end{equation*}
for any $y\in\bar{\textbf{Q}}$. This implies that  $x\in\textbf{Q}$ satisfies (\ref{equ:numnum}).

\indent{Conversely}, let $x\in\textbf{Q}$ be the solution to scheme (\ref{equ:numnum}). We need to prove that $x$ is the minimizer of $J(y)$ on $ \bar{\textbf{Q}}$.

 For any  $ y\in\partial\textbf{Q}$, $J(y)=+\infty$. We always have  $J(y)\geq J(x)$.
 Then for  any $y\in\textbf{Q}$, taking the inner product  of (\ref{equ:numnum}) with $y-x$ and using summation by parts, we get
\begin{equation}\label{equ:vp1}
  \frac{1}{\tau}\Big[\frac{f_{0}(X)}{x^{n}(1-x^{n})}(x-x^{n})\Big|y-x\Big]
    -\Big(\frac{Af_{0}(X)}{D_{h}x}\Big|D_{h}(y-x)\Big)
    +\Big[f_{0}(X)\frac{1-2x^{n}}{x^{n}(1-x^{n})}\Big|y-x\Big]=0 .
\end{equation}
After direct calculation, we see that, for  any $y\in\textbf{Q}$
  \begin{align}
  J(y)&=J(x+(y-x))\nonumber \\
  &=J(x)+\frac{1}{2\tau}\Big[\frac{f_{0}(X)}{x^{n}(1-x^{n})}(y-x)\Big|(y-x)\Big]+\frac{1}{\tau}\Big[\frac{f_{0}(X)}{x^{n}(1-x^{n})}(x-x^{n})\Big|y-x\Big]
     \nonumber\\
     &\ \ \ +\Big(Af_{0}(X)\Big|\ln\big(\frac{D_{h}x}{D_{h}y}\big)\Big)
     +\Big[f_{0}(X)\frac{1-2x^{n}}{x^{n}(1-x^{n})}\Big|y-x\Big] \nonumber  \\
  &\geq J(x), \label{equ:vp3}
  \end{align}
 where the last inequality is obtained from  (\ref{equ:vp1})   and the fact $\ln(p)\geq1-\frac{1}{p}$, for $p>0$, which leads to
 $$\Big( Af_{0}(X)\Big|\ln\big(\frac{D_{h}x}{D_{h}y}\big)\Big)\geq\Big
(Af_{0}(X)\Big|1-\frac{D_{h}y}{D_{h}x}\Big)=-\Big(\frac{Af_{0}(X)}{D_{h}x}\Big|D_{h}(y-x)\Big).$$
The proof is finished. $\hfill\Box$\\

\indent{We} define the discrete total energy $E_{N}:\textbf{Q}\rightarrow \mathbb{R}$ of \eqref{tot-energ-e} as
 $$ E_{N}(x):=\Big(Af_{0}(X)\Big|\ln\big(\frac{Af_{0}(X)}{D_{h}x}\big)\Big)+\Big[ f_{0}(X)\Big|\ln\big(x(1-x))\Big]\equiv E_{N,c}(x)- E_{N,e}(x),$$
where $ E_{N,c}(x)$ and $E_{N,e}(x)$ are both convex and their first order variations are 
\begin{equation} \label{variation} \delta_{x} E_{N,c}(x)=d_{h}\big(\frac{Af_{0}(X)}{D_{h}x}\big),\ \delta_{x}E_{N,e}(x)=-f_{0}(X)\frac{1-2x}{x(1-x)}.\end{equation}


\begin{thm}
Suppose $x^{n}=(x_{0}^{n},...,x^{n}_{N}) \in\textbf{Q}$ be the solution to scheme \eqref{equ:numnum} at time $t_{n}$. Then   the discrete  energy dissipation law holds, i.e., $$\frac{E_{N}(x^{n+1})-E_{N}(x^n)}{\Delta t } \leq - \Big[\frac{f_{0}(X)}{x^{n}(1-x^{n})}\frac{x^{n+1}-x^{n}}
{\Delta t}\Big|\frac{x^{n+1}-x^{n}}
{\Delta t}\Big], \ \ n=0,1,\cdots.$$
This is the discrete counterpart of the dissipation law in Lemma \ref{lem-dissip}.
\end{thm}
\noindent$\textbf{Proof.}$
Thanks to the convexity of $E_{N,c}^{n}$ and $E_{N,e}^{n}$, we have
\begin{align} & E_{N,c}(x^{n})-E_{N,c}(x^{n+1})\geq
  \Big[\delta_{x} E_{N,c}(x^{n+1})\Big|x^{n}-x^{n+1}\Big],\notag \\
& E_{N, e}(x^{n+1})-E_{N,e}(x^{n})\geq\Big[\delta_{x} E_{N,e}(x^n)\Big|x^{n+1}-x^{n}\Big]. \notag
\end{align}
 Then from \eqref{variation} and \eqref{equ:numnum},
\begin{subequations}
\begin{align}
E_{N}(x^{n+1})-E_{N}(x^n)&=(E_{N,c}(x^{n+1})-E_{N,e}(x^{n+1}))-(E_{N,c}(x^{n})-E_{N,e}(x^{n}))\nonumber\\
 &\leq\Big[\delta_{x}E_{N,c}(x^{n+1})-\delta_{x} E_{N,e}(x^n)\Big|x^{n+1}-x^{n}\Big] \nonumber\\
&= - \Big[\frac{f_{0}(X)}{x^{n}(1-x^{n})}\frac{x^{n+1}-x^{n}}
{\Delta t}\Big|x^{n+1}-x^{n}\Big]   \leq 0 \nonumber
\end{align}
\end{subequations}
Then the proof is completed. \hfill$\Box$

Hence the numerical scheme (\ref{equ:numnum}) for $x\in\textbf{Q}$ is uniquely solvable.  And regardless of time step, the energy decays in time: $E_{N}(x^{n+1})\leq E_{N}(x^{n})$.

Before we analyse  the convergence of damped Newton's iteration \eqref{equ:numnum1} with \eqref{Newton}, the definition of $\emph{self-concordant}$    should be involved.
\begin{defn}\label{def:self-concordant}
Let $\mathcal{G}$ be a finite-dimensional real vector space, $\mathcal{Q}$ be an open nonempty convex subset of $\mathcal{G}$, $\Lambda : \mathcal{Q}\rightarrow\mathbb{R}$ be a function, $a>0$. $\Lambda$ is called self-concordant on $\mathcal{Q}$ with the parameter value $a$, if $\Lambda \in C^3$ is a convex function on $\mathcal{Q}$, and, for all $x\in \mathcal{Q}$ and all $u\in\mathcal{G}$, the following inequality holds:
$$|D^{3}\Lambda(x)[u,u,u]|\leq 2a^{-1/2}(D^2\Lambda(x)[u,u])^{3/2}$$
($D^{k}\Lambda(x)[u_1,\cdots,u_k]$ henceforth denotes the value of the kth differential of $\Lambda$ taken at $x$ along the collection of directions $u_1,\cdots,u_k$). \cite{Y. Nesterov(1994)}
\end{defn}


\begin{thm}
Suppose $f_0(X)\in\mathcal{E}_{N}$ is the initial distribution with a positive lower bound for $X\in Q$,  then $J(y)$, defined in \eqref{DisEnergy}, is a self-concordant function  and Newton's iteration \eqref{equ:numnum1}-\eqref{Newton} is convergent in $Q$ .
\end{thm}
\noindent$\textbf{Proof.}$ Let $C_0 :=\min\limits_{X\in Q} f_0(X) >0$  and $J(y):=J_1(y)+J_2(y)+J_3(y)$ with
$$J_1(y):=\frac{1}{2\tau}\Big[\frac{f_{0}(X)}
{x^{n}(1-x^{n})}(y-x^{n})\Big|(y-x^{n})\Big],$$
$$J_2(y):=\Big
(Af_{0}(X)\Big|\ln\big(\frac{Af_{0}(X)}{D_{h}y}\big)\Big),$$
$$J_3(y):=\Big[ f_{0}(X)\frac{1-2x^{n}}{x^{n}(1-x^{n})}\Big|y\Big].$$
Since  linear and quadratic functions  have zero third derivative, $J_1(y)$ and $J_3(y)$ are  self-concordant for all $y\in Q$.  We  just need to prove $J_2(y)$ is a self-concordant function in $Q$. 

Based on the Definition \eqref{def:self-concordant}, a function $J_2 : Q\rightarrow\mathbb{R}$ is self-concordant if it is self concordant along every line in its domain,
i.e.,  $\tilde{J}_2(\xi)=J_2(y+\xi u)$ is a self-concordant function of $\xi$ for all $y\in Q$ and for all $u$  \cite{S. Boyd(2004)}.
 
 Combining with  the definition of "inner-product"  \eqref{def:inner-half},  we have
 \begin{equation}
 \tilde{J}_2(\xi )=J_2(y+\xi u)=h\sum\limits_{i=1}^{N}b_i\ln(\frac{hb_i}{y_i+\xi u_i-y_{i-1}-\xi u_{i-1}}),
  \end{equation}
and 
 \begin{equation}
\tilde{J}_2''(\xi )=h\sum\limits_{i=1}^{N}\frac{b_{i}(u_i-u_{i-1})^2}{(y_i+\xi u_i-y_{i-1}-\xi u_{i-1})^2},
 \end{equation}
and
 \begin{equation}
\tilde{J}_3'''(\xi )=-2h\sum\limits_{i=1}^{N}\frac{b_{i}(u_i-u_{i-1})^3}{(y_i+\xi u_i-y_{i-1}-\xi u_{i-1})^3},
 \end{equation}
where $h=1/N$ and $b_i=(Af_0(X))_{i-\frac{1}{2}}$, $i=1,\cdots, N$.
Then according to  the inequality: 
\begin{equation}
|\sum\limits_{i=1}^{N} w_i^3|\leq (\sum\limits_{i=1}^{N} w_i^2)^{\frac{3}{2}},\  \forall w_i\in\mathbb{R}, \nonumber
\end{equation}
 proved by Cauchy  inequality, 
 we have 
\begin{equation}
\begin{split}
  |\sum\limits_{i=1}^{N}\frac{h b_{i}(u_i-u_{i-1})^3}{(y_i+\xi u_i-y_{i-1}-\xi u_{i-1})^3}| &\leq
   \big(\sum\limits_{i=1}^{N}\frac{(hb_{i})^\frac{2}{3}(u_i-u_{i-1})^2}{(y_i+\xi u_i-y_{i-1}-\xi u_{i-1})^2}\big)^{\frac{3}{2}} \\
   &\leq \frac{1}{\sqrt{a}}\big(\sum\limits_{i=1}^{N}\frac{hb_{i}(u_i-u_{i-1})^2}{(y_i+\xi u_i-y_{i-1}-\xi u_{i-1})^2}\big)^\frac{3}{2} 
\end{split}
\end{equation}
where $a=hC_0$.   That means $J(y)$ is self-concordant for $y\in Q$.

Based on \textbf{Theorem\ 2.2.3} in \cite{Y. Nesterov(1994)},  Newton's iteration \eqref{equ:numnum1}-\eqref{Newton} is convergent in $Q$.  $\hfill\Box$

\section{Numerical Results}\label{sec:4}
\subsection{Numerical results for positive initial functions}
\begin{figure}
\captionsetup{font={scriptsize}}
\centering
\subfigure[\scriptsize $f^{1}_{0}(x)=1$]
{\includegraphics[width=6cm,height=5cm]{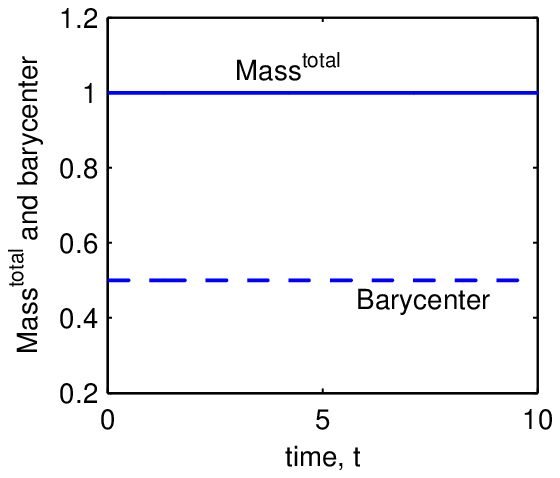}}\hspace{1in}
\subfigure[\scriptsize $f^{2}_{0}(x)=\frac{1}{5}(2+6x+\frac{\pi}{2}\sin(2\pi x))$]
{\includegraphics[width=6cm,height=5cm]{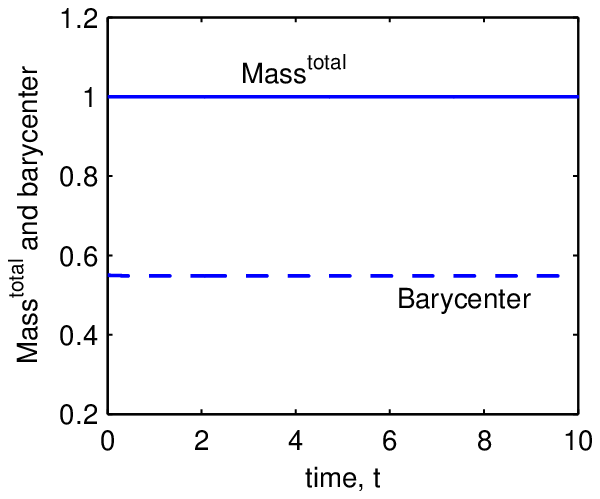}}
\caption{Total Mass ($Mass^{total}$) and Barycenter of positive initial functions over time with $h=1/1000$, $\tau=1/1000$}
\label{fig:PoMass}
\end{figure}

\begin{figure}
\captionsetup{font={scriptsize}}
\centering
\subfigure[\scriptsize $f^{1}_{0}(x)$]
{\includegraphics[width=6cm,height=5cm]{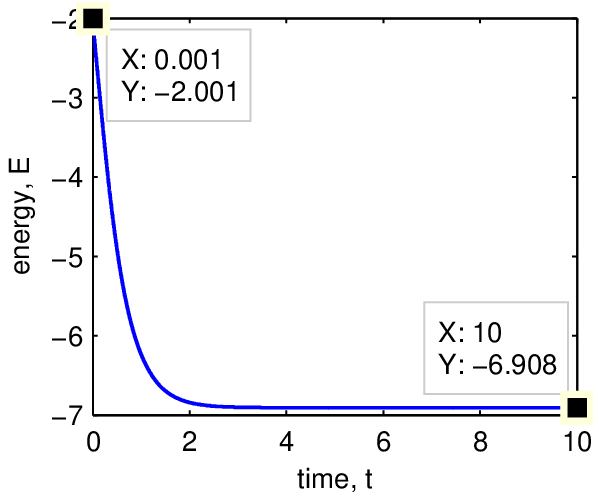}}\hspace{1in}
\subfigure[\scriptsize $f^{2}_{0}(x)$] 
{\includegraphics[width=6cm,height=5cm]{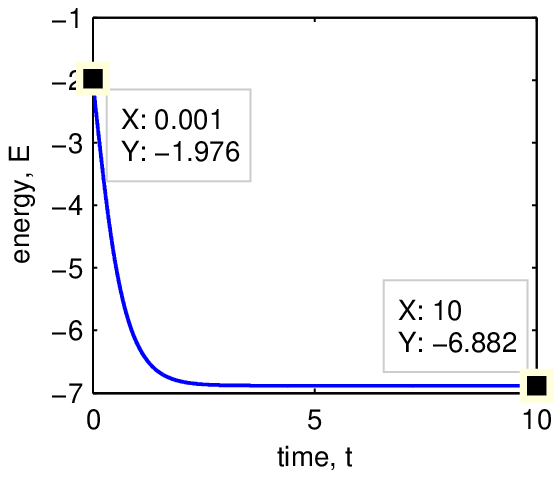}}
\caption{Energy of positive initial functions over time with $h=1/1000$, $\tau=1/1000$}
\label{fig:PoEnergy}
\end{figure}

\begin{figure}
\captionsetup{font={scriptsize}}
\centering
\subfigure[\scriptsize $t=0.002$ for $f^{1}_{0}(x)$]
{\includegraphics[width=4cm,height=3.5cm]{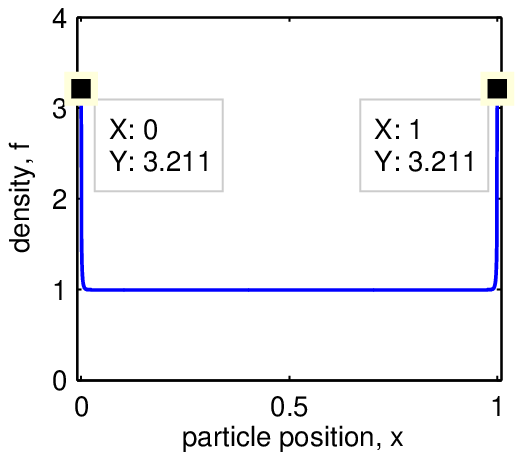}}\hspace{.3in}
\subfigure[\scriptsize $t=0.01$ for $f^{1}_{0}(x)$]
{\includegraphics[width=4cm,height=3.5cm]{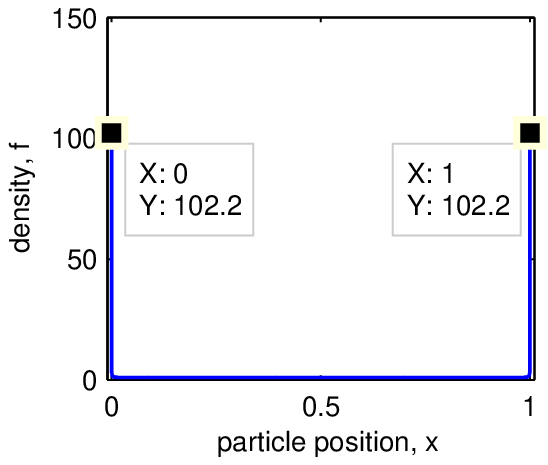}}\hspace{.3in}
\subfigure[\scriptsize  $t=10$ for $f^{1}_{0}(x)$]
{\includegraphics[width=4cm,height=3.5cm]{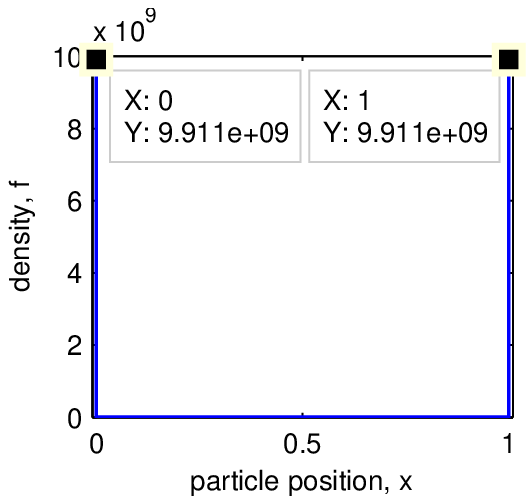}}
\caption{Density over time for  $f^{1}_{0}(x)$  with $h=1/1000$, $\tau=1/1000$}
\label{fig:density1}
\end{figure}

\begin{figure}
\captionsetup{font={scriptsize}}
\centering
\subfigure[\scriptsize  $t=0.002$ for $f^{2}_{0}(x)$] 
{\includegraphics[width=4cm,height=3.5cm]{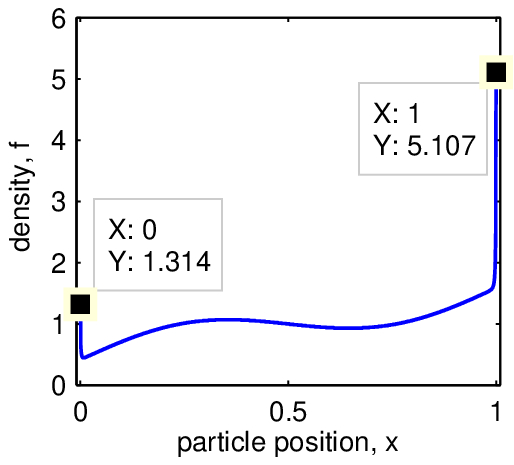}}\hspace{.3in}
\subfigure[\scriptsize  $t=0.01$ for $f^{2}_{0}(x)$]
{\includegraphics[width=4cm,height=3.5cm]{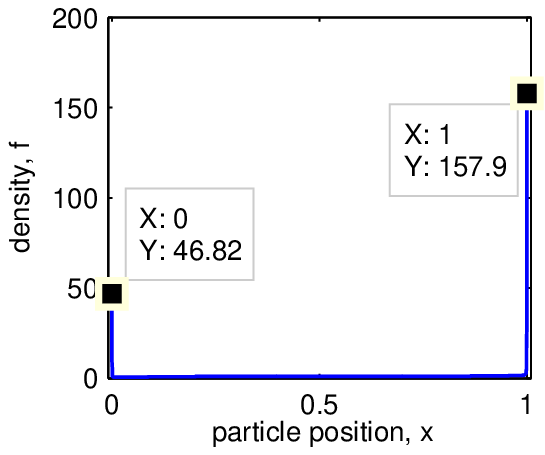}}\hspace{.3in}
\subfigure[\scriptsize  $t=10$ for $f^{2}_{0}(x)$]
{\includegraphics[width=4cm,height=3.5cm]{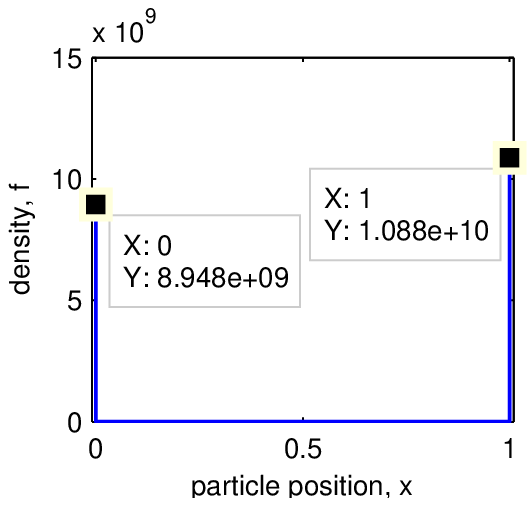}}
\caption{Density over time for   $f^{2}_{0}(x)$ with $h=1/1000$, $\tau=1/1000$}
\label{fig:density2}
\end{figure}

\begin{figure}
\captionsetup{font={scriptsize}}
\centering
\subfigure[\scriptsize $x$  for $f^{1}_{0}(x)$]
{\includegraphics[width=6cm,height=5cm]{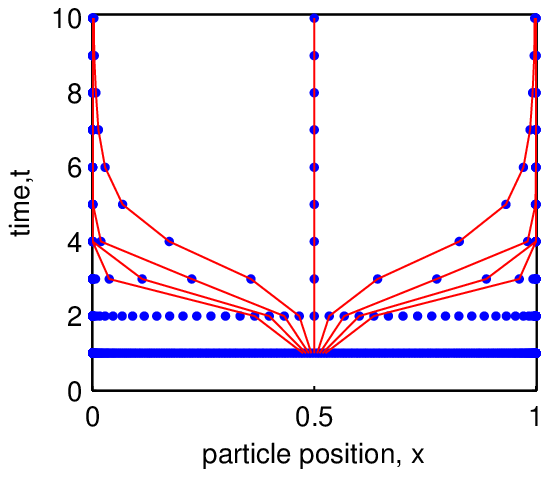}}\hspace{1in}
\subfigure[\scriptsize $x$ for $f^{2}_{0}(x)$]
{\includegraphics[width=6cm,height=5cm]{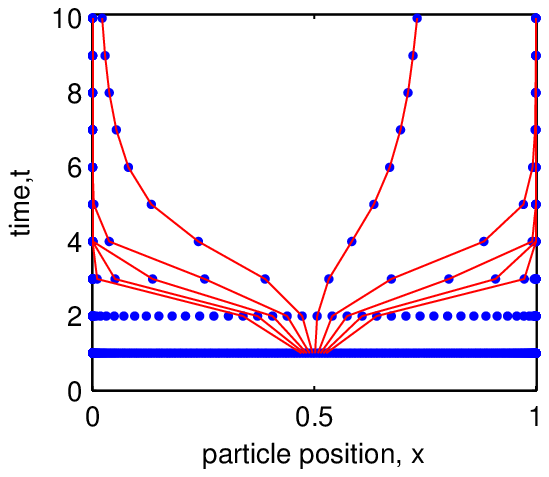}}
\caption{Particle position, $x$, over time for  $f^{1}_{0}(x)$ and $f^{2}_{0}(x)$ with $h=1/1000$, $\tau=1/1000$}
\label{fig:particle}
\end{figure}

  In this subsection, we present some numerical results for equation \eqref{eqcm}-\eqref{eqbc} with positive initial functions by Algorithm 2.1.  We take $f_{0}^{1}(x)=1$, $f_{0}^{2}(x)=\frac{1}{5}(2+6x+\frac{\pi}{2}\sin(2\pi x))$ as examples and choose the space mesh size $h=1/1000$, time step size $\tau=1/1000$  under a criterion  $\varepsilon_0=10^{-10}$. Also note that,
although the total mass of the system is equal to unity, it is not the total probability since the initial function is not  in the probability measure.  At the same time, the first moment (the mean) stands for barycenter instead of expectation. \\
\indent Fig. \ref{fig:PoMass} shows
that the total mass is unity all the time and the mean value keeps the conservation for both the positive initial functions. Fig. \ref{fig:PoEnergy} shows the total energy of the two systems decay as time evolves. The solutions of the two initial functions at time $t=0.002$,  $t=0.01$ and the steady state $t=10$ are shown in Fig. \ref{fig:density1} and Fig. \ref{fig:density2}, respectively: singularities develop at two boundaries and the heights are dependent on the mean of  initial state.  Fig. \ref{fig:particle}  shows  the motion of particles  which is influenced by the initial state. After certain time, almost all particles stay at the two boundaries, which causes $fixation\  phenomenon$. 
 This result means that we obtain the numerical complete solution, with the numerical scheme \eqref{equ:numnum} satisfying energy decay over time.  Moreover, we can  approximate the delta singularity to the scale of $1e+10$.

    Table \ref{table11} presents the total mass ($M^{total}$), barycenter (Barycenter), the  density and
the mass at the two boundary points ($f_{l}$, $f_{r}$, $M_{l}$, $M_{r}$) of the two  initial functions with different grid size ($h=1/100$, $\tau=1/100$; $h=1/1000$, $\tau=1/1000$; $h=1/10000$, $\tau=1/10000$) at time $t=10$. It shows that
the total mass keeps unity regardless of the grid size,
and  the barycenter approximates to its own initial mean at the level of the grid size.
It also shows that delta singularities at boundaries can be simulated at the level of $1e+10$ regardless of the grid size and the values are influenced by the initial expectation.
Moreover, the sum of $M_l$ and $M_r$ is approximate to unity, which verifies the development of Dirac delta functions.\\

\setcounter{equation}{0}

\begin{threeparttable}[b]
\small
\centering
\setlength{\abovecaptionskip}{10pt}
\setlength{\belowcaptionskip}{-5pt}
\captionsetup{font={scriptsize}}
\caption{\scriptsize Results for positive initial functions $f_0^1$, $f_0^2$ at time $t=10$ with different grid sizes }
\label{table11}
\begin{tabular}{@{ } l c c c c c c c}

\hline
\multicolumn{1}{l}{}
&\multicolumn{7}{c}{$f_{0}^{1}=1$}\\\hline

 $h$    &$\tau$ &$M^{total}$       & Barycenter    &$f_{l}$     &$f_{r}$   &$M_{l}$     &$M_{r}$    \\\hline
1/100 &1/100 &1.0000& 0.5000&8.2235e+09&8.2235e+09&0.4150&0.4150 \\\hline
 1/1000 &1/1000 &1.0000& 0.5000&9.9105e+09&9.9105e+09&0.4965&0.4965 \\\hline
 1/10000 &1/10000 &1.0000& 0.5000&9.9930e+09&9.9930e+09&0.4998&0.4998\\\hline

 \multicolumn{1}{l}{}
&\multicolumn{7}{c}{$f_{0}^{2}=\frac{1}{5}(2+6x+\frac{\pi}{2}\sin(2\pi x))$}\\\hline
 $h$  & $\tau$ & $M^{total}$       & Barycenter  &$f_{l}$  & $f_{r}$   & $M_{l}$ & $M_{r}$\\\hline
1/100 &1/100 &1.0000&0.5316&7.4881e+09&8.9220e+09&0.3834&0.4489 \\\hline
1/1000 &1/1000 &1.0000&0.5483&8.9477e+09&1.0879e+10&0.4475&0.5445 \\\hline
 1/10000 &1/10000 &1.0000&0.5498&8.9952e+09&1.0989e+10&0.4499&0.5496 \\\hline
\end{tabular}
\begin{tablenotes}
       \scriptsize
        \item[1]   $M^{total}$  denote by Total Mass.
        \item[2]  $M_{l}$ and $M_{r}$ are the mass  at left and right boundaries, respectively.
\end{tablenotes}
\end{threeparttable}

\subsection{Numerical results for pure drift}

In this section, we focus on $f_{0}(x)=\delta(x-x_{0})$ ($0<x_{0}<1$) and  use normal distribution $N(x_{0},\sigma^{2})$ ($\sigma=0.01$) to approximate $\delta(x-x_{0})$.
 Based on Remark \ref{reminit}, we split  the problem \eqref{eqcm}-\eqref{eqbc} into two positive initial value problems:

\begin{equation}\label{2}
   \left\{
    \begin{aligned}
     & g_t=\partial_{xx}[x(1-x)g],\ x\in (0,1),\ t>0 , \\
     & g(x,0)=10,\ x\in [0,1] , \\
     & \partial_{x}[x(1-x)g]\mid_{x=0}=0,\ \partial_{x}[x(1-x)g]\mid_{x=1}=0,\ t>0 ,
    \end{aligned}
    \right. \
\end{equation}
\begin{equation}\label{3}
   \left\{
    \begin{aligned}
     & w_t=\partial_{xx}[x(1-x)w],\ x\in (0,1),\ t>0 , \\
     & w(x,0)=10+N(x_{0},\sigma^{2}),\ x\in [0,1] , \\
     & \partial_{x}[x(1-x)w]\mid_{x=0}=0,\ \partial_{x}[x(1-x)w]\mid_{x=1}=0,\ t>0 .
    \end{aligned}
    \right. \
\end{equation}
Then we have the solution  $f=w-g$.  Because of this fact, we first obtain the numerical solutions $G(x^n, t^n)$ and $W(y^n, t^n)$ of two problems (\ref{2}) and (\ref{3}) by Algorithm 2.1, respectively, where $x^n$ and $y^n$ are the particle positions at time $t^n$. We cannot take the difference between $G(x^n, t^n)$ and $W(y^n, t^n)$ directly since $x^n$ and $y^n$ may be different. We need to get the value of $G$ at $y^n$ by the mass-conserved interpolation.
	
The details of the mass-conserved interpolation are shown as follows:

\noindent{\bf Algorithm 3.1. } {\em (Mass-conserved interpolation)
\begin{itemize}
\item Input: \ the particle positions $x=(x_0, x_1,\ldots, x_N)$ and $y=(y_0,y_1, \ldots, y_N)$; Starting point $i_s$ and ending point $i_e$ of free particles in $x$; Starting point $j_s$ and ending point $j_e$ of free particles in $y$;  Mass $m_x=(m_{x_0}, m_{x_1}, \ldots, m_{x_N})$ for each particle of $x$.\\

    \noindent{Output:} \  $m_y=(m_{y_0}, m_{y_1}, \ldots, m_{y_N})$, the re-assigned mass carried by particles $y$; $G(y)=(G_{j_s}, \ldots, G_{j_e})$, the value of  $G$ at $y$.

\item  Re-assign the mass from particles $x$ to $y$.
\begin{enumerate}
\item Define the mean mass density function $\bar{m}(s), s\in [0,1]$. Let $\Delta x_i=x_{i+1}-x_i$.
 \begin{eqnarray*}
       \bar{m}(s) & =&  \frac{m_{x_i}} {(x_{i+1}-x_{i-1})/2}, \ \mbox{for}\ s\in \left(x_i-\frac{\Delta x_{i-1}} 2, \ x_i + \frac{\Delta x_i}2\right),\  i_s <i< i_e;\\
       \bar{m}(s) & =&  \frac{m_{x_i}} {(x_{i+1}-x_i)/2}, \ \mbox{for}\ s\in \left(x_i, \ x_i + \frac{\Delta x_i}2\right), \ i=i_s; \\
       \bar{m}(s) & = & \frac{m_{x_i}} {(x_{i}-x_{i-1})/2}, \ \mbox{for}\ s\in \left(x_i-\frac{\Delta x_{i-1}} 2, \ x_i\right),\   i=i_e.
       \end{eqnarray*}
       Note that $x_{i_s} = 0$ and $x_{i_e}= 1$.
\item  Collect mass for particles at $y=(y_0, y_1, \ldots, y_N)$. Let $\Delta y_j= y_{j+1}-y_j$.\\
For free particles,
\begin{eqnarray*}
m_{y_j} &=& \int_{y_j-\frac{\Delta y_{j-1}}2}^{y_j+\frac{\Delta y_{j}}2} \bar{m}(s) ds, \ j_s < j < j_e;\\
m_{y_j} &=& \int_{y_j}^{y_j+\frac{\Delta y_{j}}2} \bar{m}(s) ds, \ j = j_s;\\
m_{y_j} &=& \int_{y_j-\frac{\Delta y_{j-1}}2}^{y_j } \bar{m}(s) ds, \  j = j_e.
\end{eqnarray*}
For particles accumulated  at left end,
\[m_{y_0}  =\cdots=m_{y_{j_s-1}}= \sum_{i=1}^{i_s-1} m_{x_i}/(i_s - 1).\]
For particles accumulated at right end,
\[m_{y_{j_e+1}}= \cdots = m_{y_N} = \sum^{N}_{i=i_e+1} m_{x_i}/(N-i_e).\]
\end{enumerate}
\item recover $G(y)=(G_{j_s}, \ldots, G_{j_e})$ from $m_{y_0}, \ldots, m_{y_N}$ by the same rules as in \eqref{numdist3}-\eqref{numdist5}.
\end{itemize}
}
%

\begin{figure}
\captionsetup{font={scriptsize}}
\centering
\subfigure[\scriptsize $x_{0}=0.4$]
{\includegraphics[width=6cm,height=5cm]{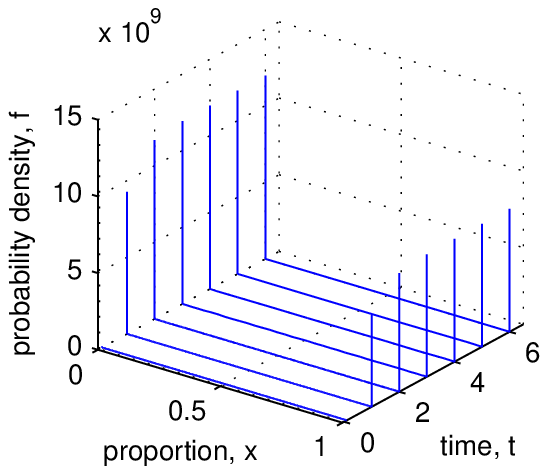}}\hspace{1in}
\subfigure[\scriptsize $x_{0}=0.7$]
{\includegraphics[width=6cm,height=5cm]{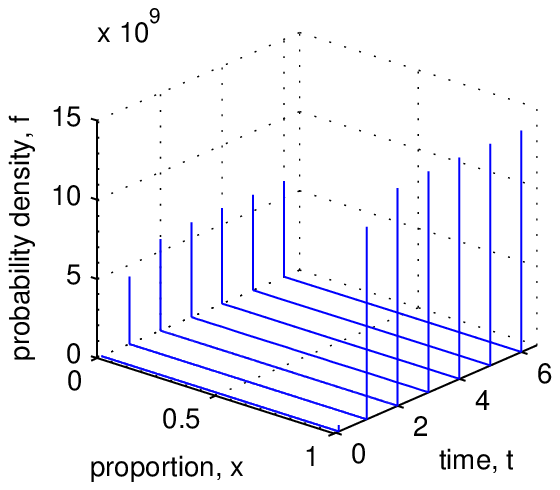}}\hspace{1in}
\caption{Distribution of probability for pure drift over time with $h=1/10000$,  $\tau=1/10000$}
\label{fig:M=0}
\end{figure}

\begin{figure}
\captionsetup{font={scriptsize}}
\centering
\subfigure[\scriptsize TP and Exp for $x_{0}=0.4$]
{\includegraphics[width=6cm,height=5cm]{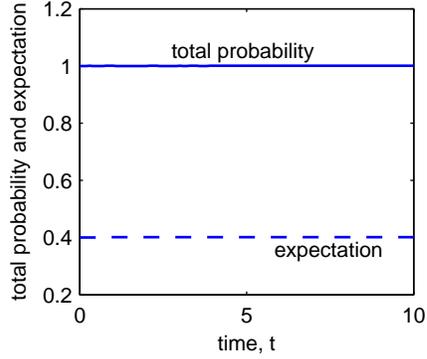}}\hspace{1in}
\subfigure[\scriptsize TP and Exp for $x_{0}=0.7$]
{\includegraphics[width=6cm,height=5cm]{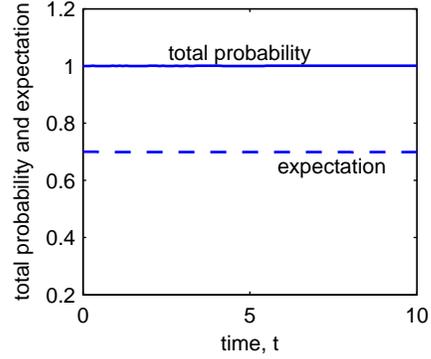}}\hspace{1in}
\caption{Total probability (TP) and expectation (Exp) for the  pure drift  as time evolves with $h=1/10000$,  $\tau=1/10000$}
\label{fig:DeltaMass}
\end{figure}

\begin{figure}
\captionsetup{font={scriptsize}}
\centering
\subfigure[\scriptsize $P_{l}$, $P_{r}$ for $x_{0}=0.4$]
{\includegraphics[width=6cm,height=5cm]{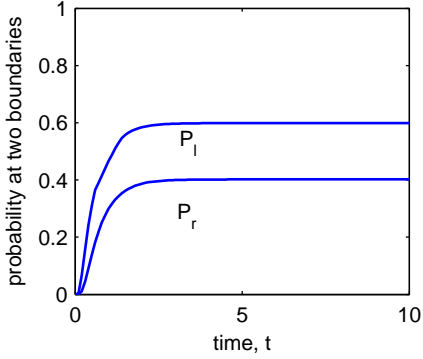}} \hspace{1in}
\subfigure[\scriptsize $P_{l}$, $P_{r}$ for $x_{0}=0.7$]
{\includegraphics[width=6cm,height=5cm]{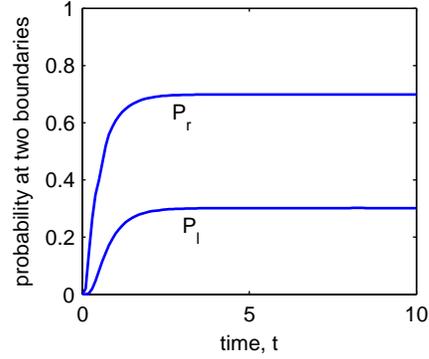}}
\caption{Probability at two boundaries  as time evolves with $h=1/10000$,  $\tau=1/10000$ ($P_{l}$ and $P_{r}$ denote the fixation probability at left  boundary and right  boundary, respectively)}
\label{fig:Mlr}
\end{figure}

Then we simulate pure drift \eqref{eqcm}-\eqref{eqbc} for $x_{0}=0.4$ and $x_{0}=0.7$  with $\varepsilon_0=10^{-10}$ and the step size $h=1/10000$,  $\tau=1/10000$. Fig.~\ref{fig:M=0} shows the evolution of distribution of probability:  the density almost vanishes in $(0,1)$, and  singularities develop at the boundary points. Moreover, the  values of singularities depend on their initial states. As shown in Fig.~\ref{fig:DeltaMass},  their total probabilities  are equal to unity and expectations keep the conservation based on their own initial expectations. This means that the numerical solution is a complete solution.
Fig.~\ref{fig:Mlr} also shows the behavior of probabilities at two boundaries  as time evolves: the value
increases to a state where the sum of both is close to unity. That causes the development of Dirac delta singularities.

Table \ref{table:two} presents the comparison of the density at two boundary points  ($f_{l}$, $f_{r}$) with scheme (3) in \cite{M. Chen(2014)}, which is a FVM scheme  with central difference method.  For $x_{0}=0.4$ and a fixed grid size $h=1/10000$, $\tau=1/10000$ with $\varepsilon_0=10^{-10}$, it shows that  $f_{l}$, $f_{r}$ obtained by scheme (3) is at the level of $1e+04$, while that scale becomes $1e+10$ by scheme \eqref{equ:numnum} in the present paper. This fact indicates that, the numerical solution obtained by scheme (\ref{equ:numnum}) is an approximation of scale $O(1/\varepsilon_0)$ to the delta singularity, with a small positive $\varepsilon_0>0$ close to the machine precision.\\


\begin{threeparttable}[b]
\small
\centering
\captionsetup{font={scriptsize}}
\caption{\scriptsize The comparision of numerical results with FVM in grid size $h=1/10000$, $\tau=1/10000$ for $x_{0}=0.4$ at $t=10$}
\label{table:two}
\begin{tabular}{@{ } l| c c| c c}
\hline
\multicolumn{1}{l|}{}
&\multicolumn{2}{c|}{ FVM }
&\multicolumn{2}{c}{Varitional Particle Scheme (\ref{equ:numnum})} \\\hline
 \ \ time &$f_{l}$    &$f_{r}$    & $f_{l}$     &$f_{r}$     \\\hline
 $\ \ t=1.0000$\ \ &\ \ 1.0039e+04\ \ &\ \ 6.0629e+03\ \ &\ \ 9.2680e+09\ \ &\ \ 5.9800e+09\ \  \\\hline
 $\ \ t=2.0000$\ \ &\ \ 1.1736e+04\ \ &\ \ 7.7362e+03\ \ &\ \ 1.1680e+10\ \ &\ \ 7.7400e+09\ \   \\\hline
 $\ \ t=3.0000$\ \ &\ \ 1.1964e+04\ \ &\ \ 7.9643e+03\ \ &\ \ 1.1930e+10\ \ &\ \ 7.7983e+09\ \   \\\hline
 $\ \ t=4.0000$\ \ &\ \ 1.1995e+04\ \ &\ \ 7.9952e+03\ \ &\ \ 1.1980e+10\ \ &\ \ 8.0200e+09\ \   \\\hline
  $\ \ t=5.0000$\ \ &\ \ 1.1999e+04\ \ &\ \ 7.9993e+03\ \ &\ \ 1.1980e+10\ \ &\ \ 8.0238e+09\ \   \\\hline
\end{tabular}
\end{threeparttable}

\subsection{Numerical results for semi-selection case}
In this part, we consider the semi-selection case where $M(x)=sx(1-x)$ ($s$ is the strength
of semi-dominant selection) in a population with the fixed size $N_e=10000$. By rescaling the time, we have the following initial-boundary value problem:
\begin{equation}\label{equ:oewM}
   \left\{
    \begin{aligned}
     & \partial_{t} f(x,t)=\frac{\partial^{2}}
     {\partial x^{2}}[x(1-x)f(x,t)]-\frac{\partial}{\partial x}[4N_{e}M(x)f(x,t)],\ x\in (0,1),\ t>0 , \\
     & f(x,0)=f_{0}(x),\ x\in [0,1] , \\
     & \{\partial_{x}[x(1-x)f(x,t)]-4N_{e}M(x)f(x,t)\}\mid_{x=0,1}=0,\ t>0 ,
    \end{aligned}
    \right. \
\end{equation}
and the corresponding energy dissipation law is given by
$$\frac{d}{dt}\left(\int_{0}^{1} f\ln\big(x(1-x)f\big)-4N_esxf dx\right)
=-\int_{0}^{1}\frac{f}{x(1-x)}|\textbf{u}|^{2}dx , $$
 where $\textbf{u}:=-\frac{\partial_{x}\left[x(1-x)f\right]}{f}+4N_{e}sx(1-x)$.
 Based on Energetic Variational Approach, Problem (\ref{equ:oewM}) is transformed into
\begin{equation}\label{equ:crpM}
   \left\{
    \begin{aligned}
     &\frac{f_{0}(X)}{x(1-x)}x_{t}=4sN_{e}f_{0}(X)-\left(\frac{\partial}{\partial X}\Big(\frac{f_{0}(X)}{\frac{\partial x}{\partial X}}\Big)
     +f_{0}(X)\frac{1-2x}{x(1-x)}\right),\ X\in(0,1),\ t>0, \\
     & x(X,0)=X,\ X\in[0,1], \\
     & x(0,t)=0,\ x(1,t)=1,\ t>0 ,
    \end{aligned}
    \right. \
\end{equation}
in the Lagrangian coordinate. Furthermore,
 the distribution of probability $\{f(x_{i}^{n}, t^{n})\}^{N}_{i=0}$ ($n>0$) can be also calculated by  \eqref{numdist3}-\eqref{numdist5}.

\begin{figure}
\captionsetup{font={scriptsize}}
\centering
\subfigure[\scriptsize $s=-0.0001$]
{\includegraphics[width=4cm,height=3.5cm]{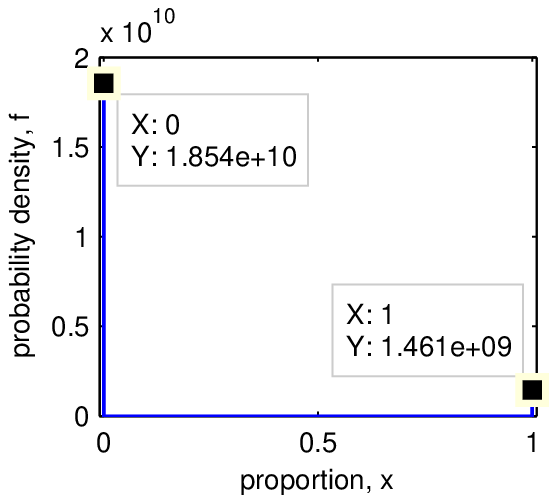}}\hspace{.3in}
\subfigure[\scriptsize $s=0.0000$]
{\includegraphics[width=4cm,height=3.5cm]{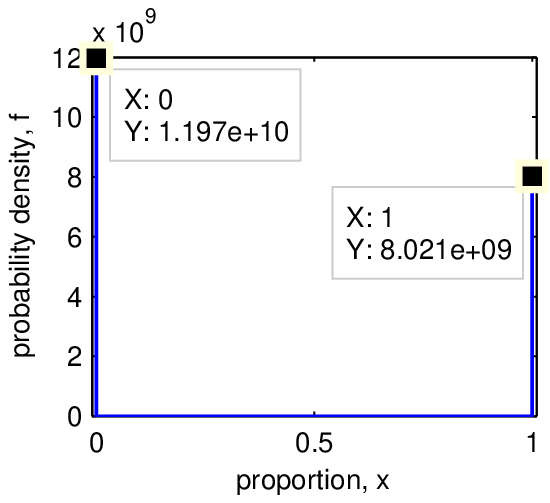}}\hspace{.3in}
\subfigure[\scriptsize $s=0.0001$]
{\includegraphics[width=4cm,height=3.5cm]{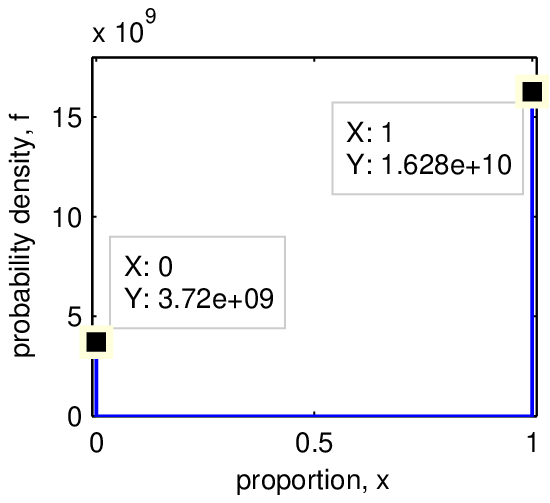}}
\caption{Distribution of probability influenced by $s$ for $x_{0}=0.4$  at  $t=10$ with $h=1/10000$ and $\tau=1/10000$}
\label{fig:den_s}
\end{figure}

\begin{figure}
\begin{minipage}[t]{0.45\linewidth}
\captionsetup{font={scriptsize}}
\centering
\includegraphics[width=6cm,height=5cm]{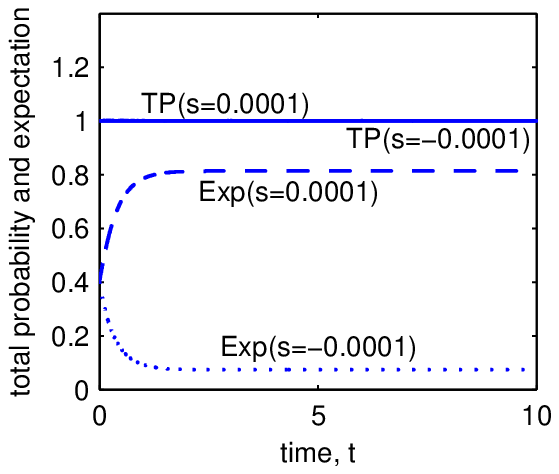}
\caption{Total probability (TP) and expectation (Exp) over time for $x_{0}=0.4$ under $s=0.0001$ and $s=-0.0001$ with $h=1/10000$ and $\tau=1/10000$}
\label{fig:S_Mass}
\end{minipage}
\begin{minipage}[t]{0.45\linewidth}
\captionsetup{font={scriptsize}}
\centering
\includegraphics[width=6cm,height=5cm]{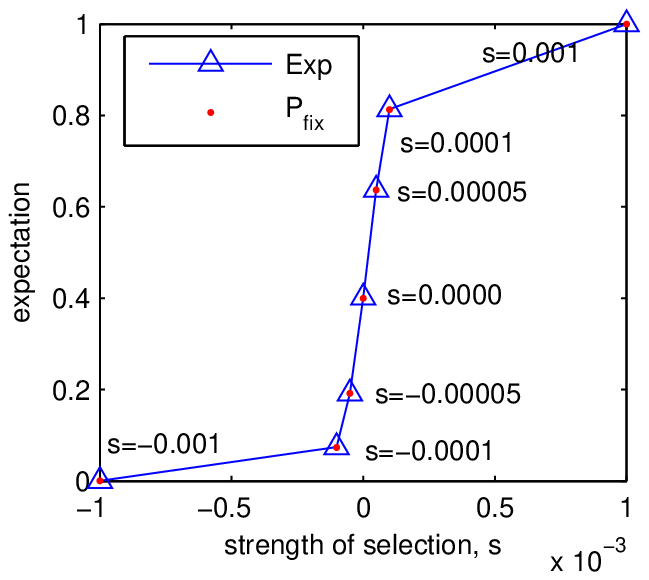}
\caption{Expectation at time $t=10$ under different $s$ for $x_{0}=0.4$  at  $t=10$ with $h=1/10000$ and $\tau=1/10000$; $P_{fix}$ is given by \eqref{P_fix} }
\label{fig:S_Exp}
\end{minipage}
\end{figure}

Fig.~\ref{fig:den_s} shows the distribution of probability of initial state $x_0=0.4$ at the steady state $t=10$ with $s=-0.0001$, $s=0.0000$, and $s=0.0001$.  It shows that semi-selection with $s=-0.0001$ prefers alleles $a$, while it is more willing to favor alleles $A$ if $s=0.0001$. Moreover, although the height of density at boundaries are influenced by $s$, they are at the scale of $1e+10$. Fig.~\ref{fig:S_Mass} implies that the total probabilities always keep normalized whatever the value of $s$ is, while the expectation does not keep conservative any more.  It means that the numerical solution  in this situation  is also a complete solution and the average is dependent on  $s$.
Fig.~\ref{fig:S_Exp} shows how the expectations are associated with
the values of $s$ when $x_0=0.4$ at time $ t=10$. It also shows that the expectation is the approximation of the probability of ultimate fixation $P_{fix}$ given by
\begin{equation}\label{P_fix}
P_{fix}(x_{0})=\frac{1-e^{-4x_{0}sN_{e}}}{1-e^{-4sN_{e}}} .
\end{equation}

\section{Conclusion and Discussion}
\label{sec:5}
In this paper, we simulate the Wright-Fisher model for pure drift and semi-selection. We first obtain the trajectory equation of the model based on EnVarA and then get the numerical scheme by the convex splitting technique. The scheme is uniquely solvable and  satisfies energy decay on a convex set where the position of particles is strictly increasing. Then we obtain the numerical  complete solutions and true probability of fixation. Moreover, at any equidistant grid, Dirac delta singularities can be measured of scale $10^{10}$ with
$\varepsilon_{0}=1e-10$ under double precision. 

 Multiple alleles at each locus among various individuals in a population, so called multiple alleles, can be considered as a high dimension problem \cite{M. Kimura(1955b), D. Waxman(2009)}. Although EnVarA can theoretically grasp the singularities on the boundary surface at a high level, it is a very challenging work to solve the constitutive relation in a high dimension. Henceforth, the numerical method based on EnVarA for the multiple alleles will be our future work.
\\
\\

 

\noindent\emph{Acknowledgments.}
It is grateful to Prof. Xinfu Chen for helpful discussions. This work is
supported in part by NSF of China under the grants 11271281. Chun Liu and Cheng Wang are partially supported by NSF grants DMS-1216938, DMS-1418689, respectively.

{

}
\end{document}